\documentclass[12pt,reqno]{amsart}
\usepackage[a4paper,margin=2.5cm,top=2.5cm,bottom=2.5cm,centering,headheight=3ex,headsep=4ex,vcentering]{geometry}

\usepackage[T1]{fontenc}
\usepackage{amsmath,amssymb,amsthm,mathtools}
\usepackage{lmodern}
\usepackage{fouriernc}
\usepackage{microtype}

\usepackage{orcidlink}
\usepackage{xcolor}

\definecolor{ao(english)}{rgb}{0.0,0.0,1.0}

\usepackage{hyperref}
\hypersetup{
  colorlinks=true,
  linkcolor=ao(english),
  citecolor=ao(english),
  urlcolor=ao(english)
}

\theoremstyle{plain}

\numberwithin{equation}{section}

\newtheorem{theorem}{Theorem}[section]
\newtheorem{lemma}[theorem]{Lemma}
\newtheorem{corollary}[theorem]{Corollary}

\newtheorem{conjecture}[theorem]{Conjecture}

\theoremstyle{definition}

\newtheorem{remark}[theorem]{Remark}

\newcommand{\PDO}{\mathrm{PDO}_{\mathrm{t}}}
\newcommand{\bx}{\mathbf{x}}

\allowdisplaybreaks

\title[Congruences modulo powers of $2$ and $3$ for Lin's function $\PDO(n)$]
{Congruences modulo powers of $2$ and $3$ for Lin's function $\PDO(n)$}



\author[M. P. Saikia]{Manjil P. Saikia\,\orcidlink{0000-0002-2997-6731}}
\address[M. P. Saikia]{Mathematical and Physical Sciences division, School of Arts \& Sciences, Ahmedabad University, Navrangpura, Ahmedabad 380009, Gujarat, India}
\email{manjil.saikia@ahduni.edu.in}

\thanks{Corresponding Author: Manjil P. Saikia (manjil.saikia@ahduni.edu.in).}

\begin{document}

\begin{abstract}
Let $\PDO(n)$ denote the total number of tagged parts over all partitions of $n$ into odd parts with designated summands. Confirming a conjecture of Singh and Barman (2025), Tang (2026) recently proved that $\PDO(3\times2^{2m-1}n)\equiv0\pmod{2^{3m+1}}$ for all $m\geq1$ and $n\geq0$. Using elementary techniques, we obtain explicit polynomial expansions for the generating functions of $\PDO(3\times2^{m+2}n)$, and we use them to determine $\PDO(3\times2^{k}n)$ modulo one further power of $2$. 

More precisely, we show that the quotient of $\PDO(3\times2^{k}n)$ by the relevant power of $2$ is odd exactly when $n$ is an odd square, or $2$, $3$ or $6$ times an odd square, depending on $k$. As a consequence, we obtain several new infinite families of congruences modulo arbitrary powers of $2$, such as
\[
\PDO\bigl(3\times4^{m}(3n+2)\bigr)\equiv0\pmod{2^{3m+4}},
\]
and
\[\PDO\bigl(3\times4^{m}(6n+5)\bigr)\equiv0\pmod{2^{3m+5}}.
\]
We also obtain an analogous improvement of a result of Li and Tang (2026) modulo powers of $3$, and an infinite family of internal congruences modulo powers of $2$.
\end{abstract}

\maketitle

\section{Introduction}\label{sec:introduction}

Andrews, Lewis, and Lovejoy~\cite{ALL} introduced partitions with designated summands, in which one occurrence of every distinct part size is distinguished (tagged) from the remaining occurrences of that size. The congruence properties of these partitions (and their generalizations) have attracted considerable attention; see, for example, \cite{HSx, BSS, SB}, etc. for some recent examples. Instead of looking at the number of such partitions, we can also look at certain statistics associated to these partitions. Lin~\cite{Lin} did exactly that and introduced the function $\PDO(n)$, which counts the total number of tagged parts over all partitions of $n$ into odd parts with designated summands, and derived the generating function
\begin{equation}\label{eq:gf}
\sum_{n\geq0}\PDO(n)q^n=q\,\frac{f_2f_3^2f_{12}^2}{f_1^2f_6}.
\end{equation}Here, and throughout this paper, we assume that $|q|<1$ and adopt the standard notations
\[
(a;q)_\infty:=\prod_{j\geq0}(1-aq^j)
\qquad\text{and}\qquad
f_k:=(q^k;q^k)_\infty\quad(k\geq1).
\]

Lin also established congruences modulo small powers of $3$ and formulated several conjectures; further arithmetic properties of $\PDO(n)$ were obtained in \cite{BK, VK, BSS}, etc. Very recently, Li and Tang~\cite{LT} proved two infinite families of congruences modulo powers of $3$ conjectured by Lin, together with two infinite families of internal congruences modulo high powers of $3$. In this paper, we are concerned with Lin's function modulo powers of $2$, for the most part. Singh and Barman~\cite[Conjecture~1.5]{SB} conjectured that, for any $k\geq0$ and $n\geq0$, we have
\begin{equation}\label{eq:conj}
\PDO\bigl(3\times2^kn\bigr)\equiv0\pmod{2^{k+2}},
\end{equation}
and they proved \eqref{eq:conj} for $0\leq k\leq6$. Very recently, Tang~\cite{T2} proved this conjecture in the following stronger form: for any $m\geq1$ and $n\geq0$, we have $\PDO(3\times2^{2m-1}n)\equiv0\pmod{2^{3m+1}}$, and this modulus is best possible. Replacing $n$ by $2n$, it follows that
\begin{equation}\label{eq:strong}
\PDO\bigl(3\times2^{2m-1}n\bigr)\equiv\PDO\bigl(3\times2^{2m}n\bigr)\equiv0\pmod{2^{3m+1}}\qquad(m\geq1,\ n\geq0).
\end{equation}

In this paper, we refine \eqref{eq:strong} by one further power of $2$, and we obtain new infinite families of congruences modulo arbitrary powers of $2$ along several arithmetic progressions. Our main tools are explicit polynomial expansions for the generating functions of $\PDO(3\times2^{m+2}n)$ for every $m\geq1$ (Theorem~\ref{thm:gf}), together with precise information on the $2$-adic orders of their coefficients (Lemma~\ref{lem:refined}).

\begin{remark}
Both Tang's proof of \eqref{eq:strong} and our approach express the relevant generating functions as $\varrho$\footnote{See Section~\ref{sec:preliminaries} for the definitions.} times a polynomial in $\xi$. Tang works with the operator $U_4=U\circ U$ (in our notation) and a modular equation of order $4$, whereas we work with the operator $U$ and the modular equation \eqref{eq:modxi} of order $2$; as illustrated in~\cite{T2}, the former modular equation can be derived from the latter. Working with $U$ gives the generating functions of $\PDO(3\times2^{m+2}n)$ for all $m\geq1$, and these consecutive expansions are what we need for the results below. The congruences \eqref{eq:strong} also follow directly from our lemmas (see Remark~\ref{rem:strong}), but they are due to Tang~\cite{T2}. None of the results stated below in this section appear in~\cite{T2}.
\end{remark}

Our first main result determines $\PDO(3\times2^{k}n)$ modulo one more power of $2$ than \eqref{eq:strong}.

\begin{theorem}\label{thm:parity}
Let $N\geq1$.
\begin{enumerate}
\item We have $\PDO(3N)\equiv0\pmod 4$, and $\PDO(3N)\equiv4\pmod 8$ if and only if $N=s^2$ or $N=3s^2$ for some odd $s$.
\item For any $m\geq1$, we have $\PDO(3\times2^{2m}N)\equiv2^{3m+1}\pmod{2^{3m+2}}$ if and only if $N=s^2$ or $N=3s^2$ for some odd $s$.
\item For any $m\geq1$, we have $\PDO(3\times2^{2m-1}N)\equiv2^{3m+1}\pmod{2^{3m+2}}$ if and only if $N=s^2$, $N=2s^2$ or $N=6s^2$ for some odd $s$.
\end{enumerate}
In cases \textup{(2)} and \textup{(3)} above, for all other $N$, we have \[\PDO(3\times2^{2m}N)\equiv\PDO(3\times2^{2m-1}N)\equiv0\pmod{2^{3m+2}}.\]\end{theorem}

Our second main result shows that, along the residue class $2$ modulo $3$, the modulus improves further.

\begin{theorem}\label{thm:mod3class}
For any $m\geq0$ and $n\geq0$, we have
\begin{align}
\PDO\bigl(3\times2^{2m}(3n+2)\bigr)&\equiv0\pmod{2^{3m+4}},\label{eq:B-3n2}\\
\PDO\bigl(3\times2^{2m}(6n+5)\bigr)&\equiv0\pmod{2^{3m+5}}.\label{eq:B-6n5}
\end{align}
\end{theorem}

As a consequence of Theorems~\ref{thm:parity} and \ref{thm:mod3class}, the modulus in \eqref{eq:strong} improves along several arithmetic progressions.

\begin{corollary}\label{cor:progressions}
For any $m\geq1$ and $n\geq0$, we have
\begin{align}
\PDO\bigl(3\times2^{2m-1}(4n+3)\bigr)&\equiv0\pmod{2^{3m+2}},\label{eq:c-4n3}\\
\PDO\bigl(3\times2^{2m-1}(8n+r)\bigr)&\equiv0\pmod{2^{3m+2}},\qquad r\in\{3,5,7\},\label{eq:c-8nr-odd}\\
\PDO\bigl(3\times2^{2m}(8n+r)\bigr)&\equiv0\pmod{2^{3m+2}},\qquad r\in\{5,7\},\label{eq:c-8nr-even}\\
\PDO\bigl(3\times2^{2m-1}(6n+5)\bigr)&\equiv\PDO\bigl(3\times2^{2m-2}(6n+5)\bigr)\equiv0\pmod{2^{3m+2}},\label{eq:c-18n15}\\
\PDO\bigl(3\times2^{2m}(3n+2)\bigr)&\equiv0\pmod{2^{3m+4}}.\label{eq:c-36n24}
\end{align}
Moreover, $\PDO(3(8n+r))\equiv0\pmod 8$ for $r\in\{5,7\}$.
\end{corollary}

\begin{corollary}\label{cor:legendre}
Let $\ell\geq5$ be a prime and let $1\leq r\leq\ell-1$ be a quadratic nonresidue modulo $\ell$.
\begin{enumerate}
\item If $\ell\equiv\pm1\pmod{12}$, then for any $n\geq0$, we have $\PDO(3(\ell n+r))\equiv0\pmod 8$ and
\[
\PDO\bigl(3\times2^{2m}(\ell n+r)\bigr)\equiv0\pmod{2^{3m+2}}\qquad(m\geq1).
\]
\item If $\ell\equiv\pm1\pmod{24}$, then for any $m\geq1$ and $n\geq0$, we have
\[
\PDO\bigl(3\times2^{2m-1}(\ell n+r)\bigr)\equiv0\pmod{2^{3m+2}}.
\]
\end{enumerate}
\end{corollary}

\begin{remark}
The restriction on $\ell$ in Corollary~\ref{cor:legendre} cannot be removed. For instance, if $\ell\equiv\pm5\pmod{12}$, then $3$ is a quadratic nonresidue modulo $\ell$, and so every nonresidue class modulo $\ell$ contains integers of the form $3s^2$ with $s$ odd, on which the first congruence fails by Theorem~\ref{thm:parity}.
\end{remark}

The same method also applies modulo powers of $3$. Li and Tang~\cite{LT} proved that $\PDO(12\times3^{m}n)\equiv0\pmod{3^{m+2}}$ for all $m\geq0$ and $n\geq0$. We show that the modulus improves outside one residue class modulo $6$.

\begin{theorem}\label{thm:three}
For any $m\geq0$ and any $N\geq0$ with $N\not\equiv1\pmod 6$, we have
\begin{equation}\label{eq:three}
\PDO\bigl(12\times3^{m}N\bigr)\equiv0\pmod{3^{m+3}}.
\end{equation}
In particular, for any $k\geq0$ and $n\geq0$, we have
\begin{equation}\label{eq:three-conj}
\PDO\bigl(4\times3^{k}(3n+2)\bigr)\equiv0\pmod{3^{k+2}}.
\end{equation}
\end{theorem}

We next consider internal congruences, that is, congruences of the form \[f(An+B)\equiv\ (\text{constant})\times f(Cn+D)\pmod M\] that relate two arithmetic subsequences of the same function. Such congruences have seen a recent flourishing, but are still relatively rarer compared to usual congruences. See, for instance the recent works \cite{CS, CT, STT, DST, LT} for some recent examples. We obtain the following infinite family of internal congruences.

\begin{theorem}\label{thm:internal}
For $k\geq4$, there exists an odd integer $\lambda_k$ such that, for any $n\geq0$, we have
\begin{equation}\label{eq:internal}
\PDO\bigl(3\times2^{k+2}n\bigr)\equiv8\lambda_k\,\PDO\bigl(3\times2^{k}n\bigr)\pmod{2^{3\lfloor k/2\rfloor+9}}.
\end{equation}
\end{theorem}

\begin{remark}
   By \eqref{eq:strong}, both sides of \eqref{eq:internal} are divisible by $2^{3\lceil k/2\rceil+4}$, so the modulus in \eqref{eq:internal} exceeds the trivial one by a factor of $2^5$ if $k$ is even and of $2^2$ if $k$ is odd. Numerical evidence suggests that much more is true (see Conjecture~\ref{conj:sharp} below).
\end{remark}

The rest of this paper is organized as follows. In Section~\ref{sec:preliminaries}, we introduce some operators and auxiliary functions that we need, and explore some of their properties. In Section~\ref{sec:U}, we describe the action of the $U$-operator on these functions, and in Section~\ref{sec:gf} we derive the generating functions of $\PDO(3\times2^{m+2}n)$. In Section~\ref{sec:proof12}, we determine the $2$-adic orders of the coefficients of these expansions. Theorems~\ref{thm:parity}, \ref{thm:mod3class} and \ref{thm:three} and Corollaries~\ref{cor:progressions} and \ref{cor:legendre} are proved in Section~\ref{sec:progressions}, and Theorem~\ref{thm:internal} is proved in Section~\ref{sec:proof3}. We conclude with a conjecture and some remarks in Section~\ref{sec:final}.

\section{Preliminaries}\label{sec:preliminaries}

For a Laurent series $\sum_{n\geq n_0}a(n)q^n$, define the operators $H$ and $U$ by
\[
H\Bigg(\sum_{n\geq n_0}a(n)q^n\Bigg):=\sum_{n\geq\lceil n_0/2\rceil}a(2n)q^{2n}
\]
and
\[
U\Bigg(\sum_{n\geq n_0}a(n)q^n\Bigg):=\sum_{n\geq\lceil n_0/2\rceil}a(2n)q^{n}.
\]
Clearly, $H(F)(q)=\frac12\bigl(F(q)+F(-q)\bigr)$, and for any two series $C$ and $F$, we have
\begin{equation}\label{eq:Uprop}
U\bigl(C(q^2)F(q)\bigr)=C(q)\,U\bigl(F(q)\bigr).
\end{equation}
We shall also use the following elementary lemma\footnote{Here, and elsewhere, we interchangeably write  $F$ and $F(q)$ for the same formal power series in $q$ depending on the context.} that is well-known in the literature.

\begin{lemma}\label{eq:conj-f}
For any $k\geq1$, we have
\[
f_k(q)f_k(-q)=
\begin{cases}
f_{2k}^3/f_{4k}, & \text{if $k$ is odd},\\[0.5ex]
f_k^2, & \text{if $k$ is even}.
\end{cases}
\]
\end{lemma}

\begin{proof}
If $k$ is even, then $f_k(-q)=f_k(q)$.

If $k$ is odd, then for every $n\geq1$ we have $(1-q^{kn})(1-(-q)^{kn})=1-q^{2kn}$ if $n$ is odd, and $(1-q^{kn})^2$ if $n$ is even. Therefore,
\[
f_k(q)f_k(-q)=\prod_{n\ \mathrm{odd}}\bigl(1-q^{2kn}\bigr)\prod_{n\ \mathrm{even}}\bigl(1-q^{kn}\bigr)^2=\frac{f_{2k}}{f_{4k}}\cdot f_{2k}^2.\qedhere
\]
\end{proof}

We need the following two $2$-dissection formulas.

\begin{lemma}\label{lem:dissections}
We have
\begin{align}
\frac{f_3}{f_1^3}&=\frac{f_4^6f_6^3}{f_2^9f_{12}^2}+3q\frac{f_4^2f_6f_{12}^2}{f_2^7},\label{eq:2diss-a}\\
\frac{f_3^3}{f_1}&=\frac{f_4^3f_6^2}{f_2^2f_{12}}+q\frac{f_{12}^3}{f_4}.\label{eq:2diss-b}
\end{align}
\end{lemma}

\begin{proof}
The first is \cite[Eq. (3.38)]{XY}, while the second is \cite[Eq. (3.75)]{XY}.
\end{proof}

We now introduce the auxiliary functions that will be used throughout this paper. Let
\begin{equation}\label{eq:aux}
\varrho:=\frac{f_2^{10}f_3^9}{f_1^{11}f_6^6},\qquad
\xi:=q\frac{f_2f_6^5}{f_1^5f_3},\qquad
\zeta:=q\frac{f_2^2f_3f_{12}^3}{f_1^3f_4f_6^2}\qquad\text{and}\qquad
\omega:=\frac{f_1^6f_6}{f_2^3f_3^2}.
\end{equation}
The functions $\varrho$ and $\xi$ are those that were also used in~\cite{LT}. All four series have integer coefficients; $\xi$ and $\zeta$ are of the form $q+O(q^2)$, and $\varrho$ and $\omega$ have constant term $1$. The generating functions considered below will be expressed as polynomials in $\xi$, multiplied by $\varrho$; the functions $\zeta$ and $\omega$ play an auxiliary role.

\begin{lemma}\label{lem:zeta}
We have
\begin{align}
1+2\zeta&=\frac{f_4f_6^7}{f_1^2f_2f_3^2f_{12}^3},\label{eq:E2}\\
1+3\zeta&=\frac{f_3f_4^3}{f_1^3f_{12}},\label{eq:E3}\\
1+4\zeta&=\frac{f_2^2f_3^4}{f_1^4f_6^2},\label{eq:E4}\\
1+6\zeta&=\frac{f_2^9f_3^2f_{12}}{f_1^6f_4^3f_6^3}.\label{eq:E6}
\end{align}
\end{lemma}

\begin{proof}
Let
\[
y:=q\frac{f_2^2f_{12}^4}{f_4^4f_6^2}.
\]
We can rewrite \eqref{eq:2diss-a} and \eqref{eq:2diss-b} as
\begin{equation}\label{eq:y-plus}
\frac{f_3}{f_1^3}=\frac{f_4^6f_6^3}{f_2^9f_{12}^2}\,(1+3y)
\qquad\text{and}\qquad
\frac{f_3^3}{f_1}=\frac{f_4^3f_6^2}{f_2^2f_{12}}\,(1+y).
\end{equation}
We now replace $q$ by $-q$ in \eqref{eq:y-plus}. Since $y$ is $q$ times a power series in $q^2$, this changes $y$ into $-y$, while the two eta quotients on the right-hand sides of \eqref{eq:y-plus} are power series in $q^2$ and remain unchanged. On the left-hand sides, Lemma~\ref{eq:conj-f} gives
\[
f_1(-q)=\frac{f_2^3}{f_1f_4}
\qquad\text{and}\qquad
f_3(-q)=\frac{f_6^3}{f_3f_{12}}.
\]
Hence \eqref{eq:y-plus} becomes
\[
\frac{f_1^3}{f_3}\cdot\frac{f_4^3f_6^3}{f_2^9f_{12}}=\frac{f_4^6f_6^3}{f_2^9f_{12}^2}\,(1-3y)
\qquad\text{and}\qquad
\frac{f_1}{f_3^3}\cdot\frac{f_4f_6^9}{f_2^3f_{12}^3}=\frac{f_4^3f_6^2}{f_2^2f_{12}}\,(1-y),
\]
that is,
\begin{equation}\label{eq:y-minus}
\frac{f_1^3}{f_3}=\frac{f_4^3}{f_{12}}\,(1-3y)
\qquad\text{and}\qquad
\frac{f_1}{f_3^3}=\frac{f_2f_4^2f_{12}^2}{f_6^7}\,(1-y).
\end{equation}
By \eqref{eq:aux} and the first identity in \eqref{eq:y-minus}, we have
\[
\zeta=y\cdot\frac{f_4^3}{f_{12}}\cdot\frac{f_3}{f_1^3}=\frac{y}{1-3y},
\]
and so
\[
1+2\zeta=\frac{1-y}{1-3y},\qquad
1+3\zeta=\frac{1}{1-3y},\qquad
1+4\zeta=\frac{1+y}{1-3y},\qquad
1+6\zeta=\frac{1+3y}{1-3y}.
\]
The first identity in \eqref{eq:y-minus} gives $\dfrac{1}{1-3y}=\dfrac{f_3f_4^3}{f_1^3f_{12}}$, which is \eqref{eq:E3}. Multiplying this by the expressions for $1-y$, $1+y$ and $1+3y$ in \eqref{eq:y-minus} and \eqref{eq:y-plus}, we obtain \eqref{eq:E2}, \eqref{eq:E4} and \eqref{eq:E6}, respectively.
\end{proof}

Let $X:=\xi(q^2)$. The relations that we shall need are collected in the following lemma.
\begin{lemma}\label{cor:relations}
We have
\begin{align}
\xi&=\zeta(1+2\zeta),\label{eq:xizeta}\\
\zeta^2&=X(1+6\zeta),\label{eq:modzeta}\\
\varrho&=\omega\,(1+8\xi)(1+9\xi),\label{eq:rho}\\
\varrho\,\xi^i&=\omega(q^2)\,\zeta^i(1+2\zeta)^i(1+3\zeta)(1+4\zeta)^2\qquad(i\geq0).\label{eq:rhoxi}
\end{align}
\end{lemma}

\begin{proof}
By \eqref{eq:aux}, the quotients $\xi/\zeta$ and $\zeta^2/X$ are equal to the right-hand sides of \eqref{eq:E2} and \eqref{eq:E6}, respectively, which gives \eqref{eq:xizeta} and \eqref{eq:modzeta}.

Next, we recall from \cite[Lemma~2.2]{LT} that
\begin{equation}\label{eq:iota}
1+8\xi=\frac{f_2^4f_3^8}{f_1^8f_6^4}
\qquad\text{and}\qquad
1+9\xi=\frac{f_2^9f_3^3}{f_1^9f_6^3}.
\end{equation} Multiplying the equations in \eqref{eq:iota} together and then with $\omega$ gives us \eqref{eq:rho}.

Finally, by \eqref{eq:xizeta}, we have $1+8\xi=(1+4\zeta)^2$ and $1+9\xi=(1+3\zeta)(1+6\zeta)$, and by \eqref{eq:aux} and \eqref{eq:E6}, we have $\omega\,(1+6\zeta)=\omega(q^2)$. Therefore, by \eqref{eq:rho} and \eqref{eq:xizeta}, we have
\begin{align*}
    \varrho\,\xi^i&= \omega \,(1+8\xi)\,(1+9\xi)\,\zeta^i(1+2\zeta)^i\\
    &=\omega\,(1+6\zeta)\cdot(1+3\zeta)(1+4\zeta)^2\cdot\zeta^i(1+2\zeta)^i\\
    &=\omega(q^2)\,\zeta^i(1+2\zeta)^i(1+3\zeta)(1+4\zeta)^2,
\end{align*}
which is \eqref{eq:rhoxi}.
\end{proof}

We use the following notation:
\[
T(q):=\sum_{\substack{s\geq1\\ s\ \mathrm{odd}}}q^{s^2}.
\]

\begin{lemma}\label{lem:mod2}
We have
\begin{align}
f_1^3f_3^3&\equiv f_4^3+qf_{12}^3\pmod 2,\label{eq:key2}\\
q^kf_{8k}^3&\equiv T(q^k)\pmod 2\qquad(k\geq1),\label{eq:jacobi2}
\end{align}
and
\begin{align}
\varrho\,\xi&\equiv T(q)+T(q^3)\pmod 2,\label{eq:rx1}\\
\varrho\,\xi^2&\equiv T(q^2)+T(q^3)+T(q^6)\pmod 2,\label{eq:rx2}\\
\varrho\,\xi^3&\equiv T(q^3)\pmod 2.\label{eq:rx3}
\end{align}
\end{lemma}

\begin{proof}
We use repeatedly that $f_k^2\equiv f_{2k}\pmod 2$. 

Multiplying the first identity in \eqref{eq:y-minus} by $f_3^4\equiv f_{12}\pmod 2$, we have
\[
f_1^3f_3^3\equiv f_4^3\,(1-3y)\equiv f_4^3+q\,\frac{f_2^2f_{12}^4}{f_4f_6^2}\equiv f_4^3+q\,\frac{f_4f_{12}^4}{f_4f_{12}}=f_4^3+qf_{12}^3\pmod 2,
\]
which is \eqref{eq:key2}. 

Next, recall Jacobi's identity \cite[Eq.~(1.7.1)]{Hir} \begin{equation}\label{jtpi}
f_1^3=\sum_{a\geq0}(-1)^a(2a+1)q^{a(a+1)/2},
\end{equation} which gives us
\[
q^kf_{8k}^3\equiv\sum_{a\geq0}q^{k(4a(a+1)+1)}=\sum_{a\geq0}q^{k(2a+1)^2}=T(q^k)\pmod 2,
\]
which is \eqref{eq:jacobi2}.

By \eqref{eq:aux}, we have
\[
\varrho\,\xi=q\,\frac{f_2^{11}f_3^8}{f_1^{16}f_6}\equiv q\,\frac{f_2^{11}f_6^4}{f_2^{8}f_6}=qf_2^3f_6^3\pmod 2.
\]
Replacing $q$ by $q^2$ in \eqref{eq:key2} and using \eqref{eq:jacobi2}, we obtain from the above
\[
\varrho\,\xi\equiv q\left(f_8^3+q^2f_{24}^3\right)\equiv T(q)+T(q^3)\pmod 2,
\]
which is \eqref{eq:rx1}. 

Further, by \eqref{eq:iota}, we have
\[
1+\xi\equiv1+9\xi=\frac{f_2^9f_3^3}{f_1^9f_6^3}\equiv\frac{f_1^9}{f_3^3}\pmod 2,
\]
and therefore, using \eqref{eq:key2} twice (the second time with $q$ replaced by $q^4$), we get
\begin{align}
\varrho\,\xi(1+\xi)&\equiv qf_1^6f_3^6\cdot\frac{f_1^9}{f_3^3}\notag\\
&=qf_1^{12}\cdot f_1^3f_3^3\equiv qf_4^3\left(f_4^3+qf_{12}^3\right)\notag\\
&\equiv qf_8^3+q^2\left(f_{16}^3+q^4f_{48}^3\right)\notag\\
&\equiv T(q)+T(q^2)+T(q^6)\pmod 2.\label{eq:rxS}
\end{align}
Subtracting \eqref{eq:rx1} from \eqref{eq:rxS} gives \eqref{eq:rx2}. Finally, $\xi^2\equiv\xi(q^2)\pmod 2$, and so
\[
\varrho\,\xi^3\equiv qf_2^3f_6^3\cdot q^2\,\frac{f_4f_{12}^5}{f_2^5f_6}=q^3\,\frac{f_4f_6^2f_{12}^5}{f_2^2}\equiv q^3f_{12}^6\equiv q^3f_{24}^3\equiv T(q^3)\pmod 2,
\]
which is \eqref{eq:rx3}.
\end{proof}

We will need the following $3$-dissections later. 

\begin{lemma}
    We have
    \begin{align}
\frac{f_2}{f_1^2}&=\frac{f_6^4f_9^6}{f_3^8f_{18}^3}+2q\,\frac{f_6^3f_9^3}{f_3^7}+4q^2\,\frac{f_6^2f_{18}^3}{f_3^6},\label{eq:3diss-a}\\
\frac{f_2^2}{f_1}&=\frac{f_6f_9^2}{f_3f_{18}}+q\,\frac{f_{18}^2}{f_9}.\label{eq:3diss-b}
\end{align}
\end{lemma}

\begin{proof}
    The first is due to Hirschhorn and Sellers \cite[Theorem 1]{HSo}, and the second can be found in \cite[Eq. (14.3.3)]{Hir}.
\end{proof}

\section{The action of the \texorpdfstring{$U$}{U}-operator}\label{sec:U}

We first describe the action of $U$ on the powers of $\zeta$.

\begin{lemma}\label{lem:Uzeta}
We have $U(\zeta)=3\xi$, $U(\zeta^2)=\xi+18\xi^2$, and, for any $i\geq3$, we have
\begin{equation}\label{eq:Uzetarec}
U\bigl(\zeta^i\bigr)=6\xi\,U\bigl(\zeta^{i-1}\bigr)+\xi\,U\bigl(\zeta^{i-2}\bigr).
\end{equation}
In particular, we have
\begin{equation}\label{eq:Uzetalist}
\begin{aligned}
U\bigl(\zeta^3\bigr)&=9\xi^2+108\xi^3, &U\bigl(\zeta^4\bigr)&=\xi^2+72\xi^3+648\xi^4,\\
U\bigl(\zeta^5\bigr)&=15\xi^3+540\xi^4+3888\xi^5, &U\bigl(\zeta^6\bigr)&=\xi^3+162\xi^4+3888\xi^5+23328\xi^6,
\end{aligned}
\end{equation}
and $U(\zeta^7)=21\xi^4+1512\xi^5+27216\xi^6+139968\xi^7$.
\end{lemma}

\begin{proof}
By \eqref{eq:modzeta}, $\zeta(q)$ is a root of $Z^2-6XZ-X=0$. Since the coefficients of this equation are series in $q^2$, replacing $q$ by $-q$ shows that $\zeta(-q)$ is also a root, and these two roots are distinct because because $\zeta(q)-\zeta(-q)=2q+O(q^3) \neq 0$. Hence $\zeta(q)+\zeta(-q)=6X$, that is, $H(\zeta)=3X$. Using \eqref{eq:modzeta}, we have
\[
    H(\zeta^2)=H(X+6X\zeta)=H(X)+6XH(\zeta)=X+18X^2.
\]
Replacing $q^2$ by $q$ gives the first two assertions.

Multiplying \eqref{eq:modzeta} by $\zeta^{i-2}$ and applying the $U$-operator, we obtain \eqref{eq:Uzetarec} from \eqref{eq:Uprop}.
\end{proof}

 Define the infinite matrix $A=(\varPhi_{i,j})_{i,j\geq1}$ as follows:
\begin{enumerate}
\item[(P1)] The first two rows of $A$ are given by
\begin{align}
&\bigl(2^4,\ 81\times2^3,\ 81\times2^6,\ 0,\ 0,\ 0,\ \ldots\bigr),\label{eq:A1}\\
&\bigl(1,\ 7\times2^5,\ 513\times2^4,\ 1539\times2^6,\ 729\times2^9,\ 0,\ 0,\ 0,\ \ldots\bigr).\label{eq:A2}
\end{align}
\item[(P2)] For any $i\geq3$ and $j\geq1$,
\begin{equation}\label{eq:A3}
\varPhi_{i,j}=10\varPhi_{i-1,j-1}+72\varPhi_{i-1,j-2}+\varPhi_{i-2,j-1}+8\varPhi_{i-2,j-2},
\end{equation}
where entries with nonpositive column indices are taken to be zero.
\end{enumerate}
By induction on $i$, we have $\varPhi_{i,j}=0$ unless $\lfloor(i+1)/2\rfloor\leq j\leq2i+1$.

\begin{lemma}\label{lem:U}
For any $i\geq1$, we have
\begin{equation}\label{eq:U}
U\bigl(\varrho\xi^i\bigr)=\varrho\sum_{j=\lfloor(i+1)/2\rfloor}^{2i+1}\varPhi_{i,j}\,\xi^j.
\end{equation}
\end{lemma}

\begin{proof}
By \eqref{eq:rhoxi}, we have
\begin{align*}
\varrho\xi&=\omega(q^2)\bigl(\zeta+13\zeta^2+62\zeta^3+128\zeta^4+96\zeta^5\bigr),\\
\varrho\xi^2&=\omega(q^2)\bigl(\zeta^2+15\zeta^3+88\zeta^4+252\zeta^5+352\zeta^6+192\zeta^7\bigr).
\end{align*}
Applying the $U$-operator and using \eqref{eq:Uprop} and Lemma~\ref{lem:Uzeta} in the above, we obtain
\begin{align*}
U(\varrho\xi)&=\omega\bigl(16\xi+920\xi^2+17352\xi^3+134784\xi^4+373248\xi^5\bigr)\\
&=\omega\,(1+8\xi)(1+9\xi)\bigl(16\xi+648\xi^2+5184\xi^3\bigr),\\
U\bigl(\varrho\xi^2\bigr)&=\omega\bigl(\xi+241\xi^2+12088\xi^3+254160\xi^4+2638656\xi^5+13436928\xi^6+26873856\xi^7\bigr)\\
&=\omega\,(1+8\xi)(1+9\xi)\bigl(\xi+224\xi^2+8208\xi^3+98496\xi^4+373248\xi^5\bigr).
\end{align*}
Using \eqref{eq:rho}, we see that, these are the cases $i=1$ and $i=2$ of \eqref{eq:U}.

Next, by \eqref{eq:xizeta} and \eqref{eq:modzeta}, $\xi=\zeta+2\zeta^2=2X+(1+12X)\zeta$. Squaring this relation and using \eqref{eq:modzeta} once more, we arrive at the modular equation
\begin{equation}\label{eq:modxi}
\xi^2=\bigl(10X+72X^2\bigr)\xi+\bigl(X+8X^2\bigr).
\end{equation}
For $i\geq3$, multiplying \eqref{eq:modxi} by $\varrho\xi^{i-2}$ and applying the $U$-operator, we find from \eqref{eq:Uprop} that
\[
U\bigl(\varrho\xi^i\bigr)=\bigl(10\xi+72\xi^2\bigr)U\bigl(\varrho\xi^{i-1}\bigr)+\bigl(\xi+8\xi^2\bigr)U\bigl(\varrho\xi^{i-2}\bigr).
\]
Together with the two initial cases and \eqref{eq:A3}, this proves \eqref{eq:U} by induction on $i$.
\end{proof}

\section{Generating functions}\label{sec:gf}

We can now find the generating function of $\PDO(3\times 2^{m+2}n)$, as given in the theorem below.
\begin{theorem}\label{thm:gf}
For any $m\geq1$, we have
\begin{equation}\label{eq:gfm}
\sum_{n\geq0}\PDO\bigl(3\times2^{m+2}n\bigr)q^n=\varrho\sum_{j=1}^{5\times2^{m-1}-1}x_{m,j}\,\xi^j,
\end{equation}
where the coefficient vectors $\bx_m=(x_{m,1},x_{m,2},\ldots)$ are defined recursively by
\begin{equation}\label{eq:x1}
\bx_1=\bigl(128\times3^3,\ 128\times3^7,\ 7168\times3^6,\ 4096\times3^8,\ 0,\ 0,\ 0,\ \ldots\bigr)
\end{equation}
and
\begin{equation}\label{eq:xrec}
\bx_{m+1}=\bx_mA.
\end{equation}
\end{theorem}

\begin{proof}
We proceed by induction on $m$. For the initial case, we start from the identity \cite[Eq.~(1.17)]{SB}, which is due to Lin~\cite{Lin}:
\begin{equation}\label{eq:L1}
\sum_{n\geq0}\PDO(6n)q^n=16q\frac{f_2^4f_3^3f_4^4}{f_1^9}=16qf_2^4f_4^4\left(\frac{f_3}{f_1^3}\right)^{3}.
\end{equation}
Substituting \eqref{eq:2diss-a} into \eqref{eq:L1}, extracting the terms whose exponents are even, and then replacing $q^2$ by $q$, we obtain
\begin{equation}\label{eq:L2}
\sum_{n\geq0}\PDO(12n)q^n=144q\frac{f_2^{18}f_3^7}{f_1^{21}f_6^2}+432q^2\frac{f_2^{10}f_3^3f_6^6}{f_1^{17}}.
\end{equation}

By \eqref{eq:aux}, we have $\omega(q^2)=\dfrac{f_2^6f_{12}}{f_4^3f_6^2}$ and using Lemma \ref{lem:zeta}, we see that
\begin{align*}
q\frac{f_2^{18}f_3^7}{f_1^{21}f_6^2}
&=\frac{f_2^6f_{12}}{f_4^3f_6^2}\cdot q\frac{f_2f_6^5}{f_1^5f_3}\cdot\left(\frac{f_3f_4^3}{f_1^3f_{12}}\right)^{2}\cdot\frac{f_2^2f_3^4}{f_1^4f_6^2}\cdot\frac{f_2^9f_3^2f_{12}}{f_1^6f_4^3f_6^3}\\
&=\omega(q^2)\,\xi\,(1+3\zeta)^2(1+4\zeta)(1+6\zeta),
\end{align*}
and
\begin{align*}
q^2\frac{f_2^{10}f_3^3f_6^6}{f_1^{17}}
&=\frac{f_2^6f_{12}}{f_4^3f_6^2}\cdot\left(q\frac{f_2f_6^5}{f_1^5f_3}\right)^{2}\cdot\frac{f_3f_4^3}{f_1^3f_{12}}\cdot\frac{f_2^2f_3^4}{f_1^4f_6^2}\\
&=\omega(q^2)\,\xi^2\,(1+3\zeta)(1+4\zeta).
\end{align*}
Since $\xi=\zeta(1+2\zeta)$ by \eqref{eq:xizeta}, substituting these two identities into \eqref{eq:L2} gives us
\begin{align*}
\sum_{n\geq0}\PDO(12n)q^n
&=144\,\omega(q^2)\,\zeta(1+2\zeta)(1+3\zeta)(1+4\zeta)\bigl\{(1+3\zeta)(1+6\zeta)+3\zeta(1+2\zeta)\bigr\}\\
&=144\,\omega(q^2)\,\zeta(1+2\zeta)(1+3\zeta)(1+4\zeta)\bigl(1+12\zeta+24\zeta^2\bigr)\\
&=144\,\omega(q^2)\bigl(\zeta+21\zeta^2+158\zeta^3+552\zeta^4+912\zeta^5+576\zeta^6\bigr).
\end{align*}
Applying the $U$-operator and using \eqref{eq:Uprop}, Lemma~\ref{lem:Uzeta} and \eqref{eq:rho}, we deduce that
\begin{align*}
\sum_{n\geq0}\PDO(24n)q^n&=144\,\omega\bigl(24\xi+2352\xi^2+71064\xi^3+943488\xi^4+5785344\xi^5+13436928\xi^6\bigr)\\
&=3456\,\omega\,(1+8\xi)(1+9\xi)\bigl(\xi+81\xi^2+1512\xi^3+7776\xi^4\bigr)\\
&=\varrho\bigl(3456\xi+279936\xi^2+5225472\xi^3+26873856\xi^4\bigr),
\end{align*}
which is \eqref{eq:gfm} for $m=1$.

Now, assume that \eqref{eq:gfm} holds for some $m\geq1$. Applying the $U$-operator to both sides and using \eqref{eq:U} and \eqref{eq:xrec}, we deduce that
\begin{align*}
\sum_{n\geq0}\PDO\bigl(3\times2^{m+3}n\bigr)q^n&=\sum_{i=1}^{5\times2^{m-1}-1}x_{m,i}\,U\bigl(\varrho\xi^i\bigr)
=\varrho\sum_{i=1}^{5\times2^{m-1}-1}x_{m,i}\sum_{j=\lfloor(i+1)/2\rfloor}^{2i+1}\varPhi_{i,j}\,\xi^j\\
&=\varrho\sum_{j=1}^{5\times2^{m}-1}\Bigg(\sum_{i=1}^{5\times2^{m-1}-1}x_{m,i}\,\varPhi_{i,j}\Bigg)\xi^j
=\varrho\sum_{j=1}^{5\times2^{m}-1}x_{m+1,j}\,\xi^j.
\end{align*}
This is \eqref{eq:gfm} with $m$ replaced by $m+1$, which completes the induction.
\end{proof}

\section{The \texorpdfstring{$2$}{2}-adic orders of the coefficients}\label{sec:proof12}

For any integer $n$, let $\nu(n)$ denote the exponent of the highest power of $2$ dividing $n$, with the convention that $\nu(0)=\infty$. We shall frequently use the elementary properties $\nu(ab)=\nu(a)+\nu(b)$ and $\nu(a_1+\cdots+a_t)\geq\min_{1\leq\ell\leq t}\nu(a_\ell)$, where equality holds in the latter if the minimum is attained by exactly one of the $a_\ell$.

We begin by establishing a lower bound for the $2$-adic orders of the entries $\varPhi_{i,j}$, which is the key technical ingredient in the proofs of our main results.

\begin{lemma}\label{lem:valPhi}
For any $i\geq1$ and $j\geq1$, we have
\begin{equation}\label{eq:valPhi}
\nu(\varPhi_{i,j})\geq2j-i.
\end{equation}
\end{lemma}

\begin{proof}
It follows immediately from \eqref{eq:A1} and \eqref{eq:A2} that
\begin{align*}
\bigl(\nu(\varPhi_{1,j})\bigr)_{j\geq1}&=(4,3,6,\infty,\infty,\ldots),\\
\bigl(\nu(\varPhi_{2,j})\bigr)_{j\geq1}&=(0,5,4,6,9,\infty,\infty,\ldots).
\end{align*}
This establishes \eqref{eq:valPhi} for $i=1$ and $i=2$.

Assume that \eqref{eq:valPhi} is valid for $1,2,\ldots,i-1$ for some $i\geq3$. Since $\nu(10)=1$, $\nu(72)=\nu(8)=3$, it follows from \eqref{eq:A3} and the inductive hypothesis that
\begin{align*}
\nu(\varPhi_{i,j})&\geq\min\bigl\{\nu(\varPhi_{i-1,j-1})+1,\ \nu(\varPhi_{i-1,j-2})+3,\ \nu(\varPhi_{i-2,j-1}),\ \nu(\varPhi_{i-2,j-2})+3\bigr\}\\
&\geq\min\bigl\{2j-i,\ 2j-i,\ 2j-i,\ 2j-i+1\bigr\}\\
&=2j-i.
\end{align*}
This completes the induction step, and therefore \eqref{eq:valPhi} holds for all $i\geq1$ and $j\geq1$.
\end{proof}

In general, the estimate \eqref{eq:valPhi} gains one power of $2$ at each application of the $U$-operator. The true behaviour of the coefficient vectors $\bx_m$ is different: their $2$-adic orders stay constant from $\bx_{2\mu+1}$ to $\bx_{2\mu+2}$ and increase by $3$ from $\bx_{2\mu+2}$ to $\bx_{2\mu+3}$. To prove this, we need the first three columns of $A$. By (P1) and (P2) of Section \ref{sec:U}, their nonzero entries, listed in increasing row order, are
\begin{equation}\label{eq:cols}
(16,\,1)^{\mathrm T},\qquad(648,\,224,\,26,\,1)^{\mathrm T}\qquad\text{and}\qquad(5184,\,8208,\,3088,\,492,\,36,\,1)^{\mathrm T},
\end{equation}
respectively. Moreover, $\varPhi_{1,j}=0$ for any $j\geq4$.

We next establish the $2$-adic orders of the coefficients $x_{m,j}$.

\begin{lemma}\label{lem:refined}
For any $\mu\geq0$, we have
\begin{align}
\nu(x_{2\mu+1,1})&=\nu(x_{2\mu+1,2})=3\mu+7, &\nu(x_{2\mu+1,j})&\geq3\mu+4+2j\quad(j\geq3),\label{eq:odd}\\
\nu(x_{2\mu+2,1})&=3\mu+7,\quad\nu(x_{2\mu+2,2})=3\mu+10, &\nu(x_{2\mu+2,j})&\geq3\mu+5+2j\quad(j\geq3).\label{eq:even}
\end{align}
\end{lemma}

\begin{proof}
By \eqref{eq:x1}, we have
\[
\bx_1=\bigl(27\times2^7,\ 2187\times2^7,\ 5103\times2^{10},\ 6561\times2^{12},\ 0,\ 0,\ \ldots\bigr),
\]
and hence $\bigl(\nu(x_{1,j})\bigr)_{j\geq1}=(7,7,10,12,\infty,\infty,\ldots)$, which implies that \eqref{eq:odd} holds for $\mu=0$.

Assume now that \eqref{eq:odd} is valid for some $\mu\geq0$. We first show that \eqref{eq:even} holds for the same $\mu$. By \eqref{eq:xrec} and \eqref{eq:cols}, we have
\begin{align*}
x_{2\mu+2,1}&=16x_{2\mu+1,1}+x_{2\mu+1,2},\\
x_{2\mu+2,2}&=648x_{2\mu+1,1}+224x_{2\mu+1,2}+26x_{2\mu+1,3}+x_{2\mu+1,4}.
\end{align*}
By the induction hypothesis, the two terms in the first line have $2$-adic orders $3\mu+11$ and $3\mu+7$, and so $\nu(x_{2\mu+2,1})=3\mu+7$. In the second line, the first term has $2$-adic order exactly $3\mu+10$, whereas the remaining three terms have $2$-adic orders $3\mu+12$, at least $3\mu+11$ and at least $3\mu+12$, respectively. Hence $\nu(x_{2\mu+2,2})=3\mu+10$. 

Now let $j\geq3$. Since $\varPhi_{1,3}=81\times2^6$ and $\varPhi_{1,j}=0$ for $j\geq4$, using the induction hypothesis together with \eqref{eq:xrec} and Lemma~\ref{lem:valPhi}, we deduce that
\begin{align*}
\nu(x_{2\mu+2,j})&\geq\min\Bigl\{\nu\bigl(x_{2\mu+1,1}\varPhi_{1,j}\bigr),\ \nu\bigl(x_{2\mu+1,2}\varPhi_{2,j}\bigr),\ \min_{i\geq3}\bigl\{\nu(x_{2\mu+1,i})+\nu(\varPhi_{i,j})\bigr\}\Bigr\}\\
&\geq\min\Bigl\{3\mu+7+2j,\ (3\mu+7)+(2j-2),\ \min_{i\geq3}\bigl\{(3\mu+4+2i)+(2j-i)\bigr\}\Bigr\}\\
&=\min\bigl\{3\mu+7+2j,\ 3\mu+5+2j,\ 3\mu+7+2j\bigr\}\\
&=3\mu+5+2j.
\end{align*}
This proves \eqref{eq:even}.

We next show that \eqref{eq:odd} holds with $\mu$ replaced by $\mu+1$. By \eqref{eq:xrec} and \eqref{eq:cols}, we have
\begin{align*}
x_{2\mu+3,1}&=16x_{2\mu+2,1}+x_{2\mu+2,2},\\
x_{2\mu+3,2}&=648x_{2\mu+2,1}+224x_{2\mu+2,2}+26x_{2\mu+2,3}+x_{2\mu+2,4}.
\end{align*}
By \eqref{eq:even}, the two terms in the first line have $2$-adic orders $3\mu+11$ and $3\mu+10$, and so $\nu(x_{2\mu+3,1})=3\mu+10$. In the second line, the first term has $2$-adic order exactly $3\mu+10$, whereas the remaining three terms have $2$-adic orders $3\mu+15$, at least $3\mu+12$ and at least $3\mu+13$, respectively. Hence $\nu(x_{2\mu+3,2})=3\mu+10$. For the third coordinate, by \eqref{eq:cols} and \eqref{eq:even}, we have
\begin{align*}
\nu(x_{2\mu+3,3})&\geq\min_{1\leq i\leq6}\bigl\{\nu(x_{2\mu+2,i})+\nu(\varPhi_{i,3})\bigr\}\\
&\geq\min\bigl\{3\mu+13,\ 3\mu+14,\ 3\mu+15,\ 3\mu+15,\ 3\mu+17,\ 3\mu+17\bigr\}\\
&=3\mu+13.
\end{align*}

Finally, let $j\geq4$. Since $\varPhi_{1,j}=0$, using \eqref{eq:even} together with \eqref{eq:xrec} and Lemma~\ref{lem:valPhi}, we deduce that
\begin{align*}
\nu(x_{2\mu+3,j})&\geq\min\Bigl\{\nu\bigl(x_{2\mu+2,2}\varPhi_{2,j}\bigr),\ \min_{i\geq3}\bigl\{\nu(x_{2\mu+2,i})+\nu(\varPhi_{i,j})\bigr\}\Bigr\}\\
&\geq\min\Bigl\{(3\mu+10)+(2j-2),\ \min_{i\geq3}\bigl\{(3\mu+5+2i)+(2j-i)\bigr\}\Bigr\}\\
&=3\mu+8+2j.
\end{align*}
Therefore, $\nu(x_{2\mu+3,j})\geq3(\mu+1)+4+2j$ for any $j\geq3$. This completes the induction step, and hence the proof.
\end{proof}

\begin{remark}\label{rem:strong}
Lemma~\ref{lem:refined} gives a short proof of Tang's congruences \eqref{eq:strong}. Indeed, let $\mu\geq0$. By Lemma~\ref{lem:refined}, we have $\nu(x_{2\mu+1,j})\geq3\mu+7$ and $\nu(x_{2\mu+2,j})\geq3\mu+7$ for any $j\geq1$. Since $\varrho$ and $\xi$ have integer coefficients, so does $\varrho\xi^j$ for every $j\geq1$, and so \eqref{eq:gfm} gives
\[
\sum_{n\geq0}\PDO\bigl(3\times2^{2\mu+3}n\bigr)q^n\equiv\sum_{n\geq0}\PDO\bigl(3\times2^{2\mu+4}n\bigr)q^n\equiv0\pmod{2^{3\mu+7}},
\]
which is \eqref{eq:strong} for $m=\mu+2\geq2$. For $m=1$, \eqref{eq:strong} follows from \eqref{eq:L1} and \eqref{eq:L2}. 
\end{remark}

\section{Proofs of Theorems~\ref{thm:parity}, \ref{thm:mod3class} and \ref{thm:three}}\label{sec:progressions}

In this section we prove Theorems~\ref{thm:parity}, \ref{thm:mod3class} and \ref{thm:three}, together with Corollaries~\ref{cor:progressions} and \ref{cor:legendre}. The idea of the proofs is the same as in \cite[Section~5]{DST}.

We first need some preliminary results.
\begin{lemma}\label{lem:mod16}
For any $N\equiv2\pmod 3$, the coefficient of $q^N$ in $\varrho\,\xi$ is divisible by $16$.
\end{lemma}

\begin{proof}
By \eqref{eq:aux}, we have
\[
\varrho\,\xi=q\,\frac{f_2^{11}f_3^8}{f_1^{16}f_6}=q\,f_2^3\cdot\frac{f_3^8}{f_6}\cdot\left(\frac{f_2}{f_1^2}\right)^{8}=:\sum_{n\geq0}a_nq^{n+1},
\]
so the coefficient of $q^N$ in $\varrho\,\xi$ is $a_{N-1}$, and we have to show that
\[
a_{3n+1}\equiv0\pmod{16}\qquad(n\geq0).
\]

Let ${A}$, ${B}$ and ${C}$ denote the three eta quotients on the right-hand side of \eqref{eq:3diss-a}, so that
\[
\frac{f_2}{f_1^2}={A}+2{R},\qquad\text{where}\quad{R}:=q\,{B}+2q^2\,{C}.
\]
Here ${A}$, ${B}$ and ${C}$ are power series in $q^3$ with integer coefficients. By the binomial theorem, from above, we have
\[
\left(\frac{f_2}{f_1^2}\right)^{8}={A}^8+\sum_{i=1}^{8}\binom{8}{i}2^i\,{A}^{8-i}{R}^i.
\]
Since $\binom{8}{1}2=16$, $\binom{8}{2}2^2=16\times7$, $\binom{8}{3}2^3=16\times28$ and $2^i\equiv0\pmod{16}$ for $i\geq4$, every term of the sum is divisible by $16$. Hence
\begin{equation}\label{nn1}
    \sum_{n\geq0}a_nq^n\equiv f_2^3\cdot\frac{f_3^8}{f_6}\,{A}^8\pmod{16}.
\end{equation}

Now write
\begin{equation}\label{nn2}
  \frac{f_3^8}{f_6}\,{A}^8=:\sum_{k\geq0}b_kq^{3k},  
\end{equation}
which is possible since the left-hand side is a power series in $q^3$ with integer coefficients. By \eqref{jtpi}, we have
\[
f_2^3=\sum_{j\geq0}(-1)^j(2j+1)\,q^{j(j+1)}.
\]
Comparing the coefficients of $q^n$ using \eqref{nn1} and \eqref{nn2}, we obtain
\[
a_n\equiv\sum_{\substack{j,k\geq0\\ 3k+j(j+1)=n}}(-1)^j(2j+1)\,b_k\pmod{16}.
\]
Finally, $j(j+1)\equiv0,2,0\pmod 3$ according as $j\equiv0,1,2\pmod 3$. Thus $3k+j(j+1)\not\equiv1\pmod 3$ for all $j,k\geq0$, and if $n\equiv1\pmod 3$, the sum is empty. Therefore, $a_{3n+1}\equiv0\pmod{16}$ for any $n\geq0$, as required.
\end{proof}

The next lemma deals with the three initial cases which are not covered by Theorem~\ref{thm:gf}.

\begin{lemma}\label{lem:base}
We have
\begin{align}
\sum_{n\geq0}\PDO(3n)q^n&=4q\,\frac{f_2f_4^2f_6^3}{f_1^4}\equiv4\bigl(T(q)+T(q^3)\bigr)\pmod 8,\label{eq:B0}\\
\sum_{n\geq0}\PDO(9n+6)q^{3n+2}&\equiv16q^2f_{48}\pmod{32},\label{eq:B0b}\\
\sum_{n\geq0}\PDO(6n)q^n&\equiv16\bigl(T(q)+T(q^2)+T(q^6)\bigr)\pmod{32},\label{eq:B1}\\
\sum_{n\geq0}\PDO(12n)q^n&\equiv144\,\varrho\left(\xi+8\xi^2+8\xi^3\right)\pmod{2^8}.\label{eq:B2}
\end{align}
\end{lemma}

\begin{proof}
By \eqref{eq:gf}, we have $\sum\limits_{n\geq0}\PDO(n)q^n=q\,\dfrac{f_3^2f_{12}^2}{f_6}\cdot\dfrac{f_2}{f_1^2}$. Substituting \eqref{eq:3diss-a}, extracting the terms in which the exponent of $q$ is a multiple of $3$, and replacing $q^3$ by $q$, we obtain the equality in \eqref{eq:B0}. 

Moreover,
\[
q\,\frac{f_2f_4^2f_6^3}{f_1^4}\equiv q\,\frac{f_2^5f_6^3}{f_2^2}=qf_2^3f_6^3\pmod 2,
\]
and we have seen in the proof of Lemma~\ref{lem:mod2} that $qf_2^3f_6^3\equiv T(q)+T(q^3)\pmod 2$. This proves \eqref{eq:B0}.

For \eqref{eq:B0b}, let ${A}$, $2q{B}$, $4q^2{C}$ denote the three terms on the right-hand side of \eqref{eq:3diss-a}, and let ${P}+q^2{Q}$ be the right-hand side of \eqref{eq:3diss-b} with $q$ replaced by $q^2$, so that ${A},{B},{C},{P},{Q}$ are power series in $q^3$. By \eqref{eq:B0}, we have
\[
\sum_{n\geq0}\PDO(3n)q^n=4qf_6^3\cdot\frac{f_4^2}{f_2}\cdot\left(\frac{f_2}{f_1^2}\right)^{2}
=4qf_6^3\left({P}+q^2{Q}\right)\left({A}+2q{B}+4q^2{C}\right)^{2}.
\]
Extracting the terms in which the exponent of $q$ is congruent to $2$ modulo $3$, we obtain
\begin{equation}\label{eq:B0c}
\sum_{n\geq0}\PDO(9n+6)q^{3n+2}=16q^2f_6^3\left({A}{B}{P}+q^3{B}^2{Q}+2q^3{A}{C}{Q}+4q^3{C}^2{P}\right).
\end{equation}

Now, modulo $2$, we have
\[
{A}\equiv1,\qquad{B}\equiv\frac{f_9^3}{f_3},\qquad{P}=\frac{f_{12}f_{18}^2}{f_6f_{36}}\equiv f_6,\qquad{Q}=\frac{f_{36}^2}{f_{18}}\equiv f_{18}^3,
\]
and, dividing \eqref{eq:key2} by $f_1^4\equiv f_4\pmod 2$ and replacing $q$ by $q^3$,
\[
\frac{f_9^3}{f_3}\equiv f_{12}^2+q^3\,\frac{f_{36}^3}{f_{12}}\equiv f_{12}^2+q^3\,\frac{f_{18}^6}{f_6^2}\pmod 2.
\]
Therefore, since ${B}^2{Q}\equiv\dfrac{f_{18}^3}{f_6}\cdot f_{18}^3\pmod 2$, we have
\[
{A}{B}{P}+q^3{B}^2{Q}\equiv f_6\left(f_{12}^2+q^3\,\frac{f_{18}^6}{f_6^2}\right)+q^3\,\frac{f_{18}^6}{f_6}\equiv f_6f_{12}^2\pmod 2.
\]
Since $f_6^4f_{12}^2\equiv f_{48}\pmod 2$, \eqref{eq:B0b} follows from \eqref{eq:B0c}.

Next, by \eqref{eq:L1} and \eqref{eq:rxS}, we have
\[
\frac{1}{16}\sum_{n\geq0}\PDO(6n)q^n=q\,\frac{f_2^4f_3^3f_4^4}{f_1^9}\equiv qf_1^{15}f_3^3\equiv T(q)+T(q^2)+T(q^6)\pmod 2,
\]
which is \eqref{eq:B1}.

Finally, in the proof of Theorem~\ref{thm:gf} we have seen that
\[
\sum_{n\geq0}\PDO(12n)q^n=\omega(q^2)\,\xi\,(1+3\zeta)(1+4\zeta)\Bigl(144(1+3\zeta)(1+6\zeta)+432\,\xi\Bigr).
\]
Since $(1+3\zeta)(1+6\zeta)=1+9\xi$, the case $i=1$ of \eqref{eq:rhoxi} gives
\begin{equation}\label{eq:L2exact}
(1+4\zeta)\sum_{n\geq0}\PDO(12n)q^n=144\,\varrho\,\xi\,(1+12\xi).
\end{equation}
By \eqref{eq:xizeta}, we have $\zeta=\xi-2\zeta^2\equiv\xi-2\xi^2\pmod 4$, and so
\[
\frac{1+12\xi}{1+4\zeta}\equiv(1+12\xi)(1-4\zeta)\equiv1+12\xi-4\xi+8\xi^2\equiv1+8\xi+8\xi^2\pmod{16}.
\]
Since $144=2^4\times9$, \eqref{eq:B2} follows from \eqref{eq:L2exact}.
\end{proof}

\begin{lemma}\label{lem:reduce}
For any $m\geq2$, there exist integers $a_m$, $b_m$, $c_m$, $d_m$ with
\[
\nu(a_m)=\nu(b_m)=\nu(c_m)=3m+1\qquad\text{and}\qquad\nu(d_m)=3m+4
\]
such that
\begin{align}
\sum_{n\geq0}\PDO\bigl(3\times2^{2m-1}n\bigr)q^n&\equiv\varrho\left(a_m\,\xi+b_m\,\xi^2\right)\pmod{2^{3m+4}},\label{eq:red-odd}\\
\sum_{n\geq0}\PDO\bigl(3\times2^{2m}n\bigr)q^n&\equiv\varrho\left(c_m\,\xi+d_m\,\xi^2\right)\pmod{2^{3m+5}}.\label{eq:red-even}
\end{align}
\end{lemma}

\begin{proof}
We take $\mu=m-2$ in Lemma~\ref{lem:refined}. By \eqref{eq:odd}, we have $\nu(x_{2m-3,1})=\nu(x_{2m-3,2})=3m+1$ and $\nu(x_{2m-3,j})\geq3m-2+2j\geq3m+4$ for $j\geq3$, and so \eqref{eq:red-odd} follows from \eqref{eq:gfm} with $(a_m,b_m)=(x_{2m-3,1},x_{2m-3,2})$. 

Similarly, by \eqref{eq:even}, we have $\nu(x_{2m-2,1})=3m+1$, $\nu(x_{2m-2,2})=3m+4$ and $\nu(x_{2m-2,j})\geq3m-1+2j\geq3m+5$ for $j\geq3$, and so \eqref{eq:red-even} follows from \eqref{eq:gfm} with $(c_m,d_m)=(x_{2m-2,1},x_{2m-2,2})$.
\end{proof}

We can now complete the proofs.

\begin{proof}[Proof of Theorem~\ref{thm:parity}]
Part (1) is just \eqref{eq:B0}. 

For part (2) with $m=1$, \eqref{eq:B2} and \eqref{eq:rx1} give
\[
\sum_{n\geq0}\PDO(12n)q^n\equiv16\,\varrho\,\xi\equiv16\bigl(T(q)+T(q^3)\bigr)\pmod{32},
\]
and for $m\geq2$, \eqref{eq:red-even} and \eqref{eq:rx1} give
\[
\sum_{n\geq0}\PDO\bigl(3\times2^{2m}n\bigr)q^n\equiv2^{3m+1}\varrho\,\xi\equiv2^{3m+1}\bigl(T(q)+T(q^3)\bigr)\pmod{2^{3m+2}}.
\]

For part (3), the case $m=1$ is \eqref{eq:B1}, and for $m\geq2$, since $a_m/2^{3m+1}$ and $b_m/2^{3m+1}$ are both odd, \eqref{eq:red-odd} and \eqref{eq:rxS} give
\[
\sum_{n\geq0}\PDO\bigl(3\times2^{2m-1}n\bigr)q^n\equiv2^{3m+1}\varrho\,\xi(1+\xi)\equiv2^{3m+1}\bigl(T(q)+T(q^2)+T(q^6)\bigr)\pmod{2^{3m+2}}.
\]
Since the sets $\{s^2\}$, $\{2s^2\}$, $\{3s^2\}$, $\{6s^2\}$ ($s$ odd) are pairwise disjoint, the result follows.
\end{proof}

\begin{proof}[Proof of Theorem~\ref{thm:mod3class}]
For $m=0$, both congruences follow from \eqref{eq:B0b}, since every exponent of $q$ occurring in $q^2f_{48}$ is even.

Let $N\equiv2\pmod 3$, and let $\alpha(N)$, $\beta(N)$ and $\gamma(N)$ denote the coefficients of $q^N$ in $\varrho\,\xi$, $\varrho\,\xi^2$ and $\varrho\,\xi^3$, respectively. By Lemma~\ref{lem:mod16}, we have $\alpha(N)\equiv0\pmod{16}$. If, in addition, $N$ is odd, then $N$ is neither of the form $2s^2$ or $6s^2$, nor of the form $3s^2$, and so by \eqref{eq:rx2} and \eqref{eq:rx3}, both $\beta(N)$ and $\gamma(N)$ are even.

For $m=1$, \eqref{eq:B2} gives
\[
\PDO(12N)\equiv2^4\times9\,\alpha(N)+2^7\times9\,\bigl(\beta(N)+\gamma(N)\bigr)\pmod{2^8},
\]
which is divisible by $2^7$, and by $2^8$ if $N$ is odd. For $m\geq2$, \eqref{eq:red-even} gives
\[
\PDO\bigl(3\times2^{2m}N\bigr)\equiv c_m\,\alpha(N)+d_m\,\beta(N)\pmod{2^{3m+5}},
\]
where $\nu(c_m\alpha(N))\geq3m+5$ and $\nu(d_m\beta(N))\geq3m+4$, with $\nu(d_m\beta(N))\geq3m+5$ if $N$ is odd. This completes the proof.
\end{proof}

\begin{proof}[Proof of Corollary~\ref{cor:progressions}]
An odd square is congruent to $1$ modulo $8$. Hence, if $s$ is odd, then $s^2\equiv1$ and $3s^2\equiv3\pmod 8$, while $2s^2$ and $6s^2$ are even. Thus $4n+3$ and $8n+r$ with $r\in\{3,5,7\}$ are never of the form $s^2$, $2s^2$ or $6s^2$, and $8n+r$ with $r\in\{5,7\}$ is never of the form $s^2$ or $3s^2$. This gives \eqref{eq:c-4n3}--\eqref{eq:c-8nr-even} and the last assertion by Theorem~\ref{thm:parity}.

Similarly, $6n+5$ is odd and congruent to $2$ modulo $3$, and so it is not of the form $s^2$, $2s^2$ or $6s^2$; this gives the first congruence in \eqref{eq:c-18n15}. The second congruence in \eqref{eq:c-18n15} is \eqref{eq:B-6n5} with $m$ replaced by $m-1$, and \eqref{eq:c-36n24} is \eqref{eq:B-3n2}.
\end{proof}

\begin{proof}[Proof of Corollary~\ref{cor:legendre}]
If $\ell\equiv\pm1\pmod{12}$, then $3$ is a quadratic residue modulo $\ell$, and so $s^2$ and $3s^2$ are quadratic residues or multiples of $\ell$ for every $s$. If $\ell\equiv\pm1\pmod{24}$, then $2$ and $3$ are both quadratic residues modulo $\ell$, and the same holds for $s^2$, $2s^2$ and $6s^2$. In either case $\ell n+r$ is not of the relevant form, and the result follows from Theorem~\ref{thm:parity}.
\end{proof}

For the proof of Theorem~\ref{thm:three} we now focus on the generating functions of Li and Tang \cite{LT} modulo powers of $3$. Let $\vartheta:=\dfrac{f_2^8f_3^5}{f_1^7f_6^4}$ and let $\nu_3$ denote the $3$-adic valuation.

\begin{proof}[Proof of Theorem~\ref{thm:three}]
By \cite[Theorem~3.2]{LT}, for any $m\geq0$, we have
\[
\sum_{n\geq0}\PDO\bigl(12\times3^{m}n\bigr)q^n=\vartheta\sum_{j\geq1}y_{m,j}\,\xi^j
\]
for certain integers $y_{m,j}$, and by \cite[Lemma~4.4]{LT}, we have $\nu_3(y_{m,1})=m+2$ and $\nu_3(y_{m,j})\geq m+3$ for $j\geq2$. Hence,
\[
\sum_{n\geq0}\PDO\bigl(12\times3^{m}n\bigr)q^n\equiv y_{m,1}\,\vartheta\,\xi\pmod{3^{m+3}}.
\]
Since $f_k^3\equiv f_{3k}\pmod 3$, we have
\[
\vartheta\,\xi=q\,\frac{f_2^9f_3^4f_6}{f_1^{12}}\equiv q\,\frac{f_6^3f_3^4f_6}{f_3^4}=qf_6^4\pmod 3,
\]
in which every exponent of $q$ is congruent to $1$ modulo $6$. This proves \eqref{eq:three}. 

For $k\geq1$, \eqref{eq:three-conj} is the case $m=k-1$, $N=3n+2$ of \eqref{eq:three}, and for $k=0$ it is the congruence $\PDO(12n+8)\equiv0\pmod 9$ of Lin \cite{Lin}.
\end{proof}

\begin{remark}
The proof shows more: for any $m\geq0$ and $n\geq0$, we have
\[
\PDO\bigl(12\times3^{m}(6n+1)\bigr)\equiv y_{m,1}\,\tau(n)\pmod{3^{m+3}},\qquad\text{where}\quad\sum_{n\geq0}\tau(n)q^n=f_1^4 .
\]
\end{remark}

\section{Proof of Theorem~\ref{thm:internal}}\label{sec:proof3}

Following the determinant method of Atkin and O'Brien~\cite{AO}, as in~\cite{LT}, we give a proof of Theorem \ref{thm:internal} in this section. We refer the reader to \cite{LT} for more details. In view of the two-step pattern in Lemma~\ref{lem:refined}, we compare $\bx_{m+2}$ with $\bx_m$. For any $m\geq1$, $s\geq1$ and $t\geq1$, let us define
\begin{equation}\label{eq:Gamma}
\varGamma^m_{s,t}:=x_{m+2,s}\,x_{m,t}-x_{m,s}\,x_{m+2,t}.
\end{equation}
Notice that $\varGamma^m_{s,s}=0$ and $\varGamma^m_{s,t}=-\varGamma^m_{t,s}$.

\begin{lemma}\label{lem:Gamma}
For any $m\geq2$, $s\geq1$ and $t\geq1$,
\begin{equation}\label{eq:valGamma}
\nu\bigl(\varGamma^m_{s,t}\bigr)\geq3m+10+2s+2t.
\end{equation}
\end{lemma}

\begin{proof}
It suffices to consider $s>t$. We prove \eqref{eq:valGamma} by induction on $m$.

For the base case $m=2$, the vectors $\bx_2$ and $\bx_4$ have $9$ and $39$ nonzero coordinates, respectively, and they are computed from \eqref{eq:x1} and \eqref{eq:xrec}. If $t\geq10$, then $\varGamma^2_{s,t}=0$, and if $s\geq10>t$, then $\varGamma^2_{s,t}=x_{4,s}x_{2,t}$. A direct computation shows that $\nu(\varGamma^2_{s,t})\geq16+2s+2t$ for all $s>t$ with $t\leq9$ and $s\leq39$, with equality for $(s,t)=(3,1)$.

Assume that \eqref{eq:valGamma} holds for some $m\geq2$. Using \eqref{eq:xrec} for both $\bx_{m+3}=\bx_{m+2}A$ and $\bx_{m+1}=\bx_mA$, we obtain
\begin{equation}\label{eq:Gammarec}
\varGamma^{m+1}_{s,t}=\sum_{k\geq1}\sum_{\ell\geq1}\varGamma^m_{k,\ell}\,\varPhi_{k,s}\varPhi_{\ell,t}
=\sum_{k>\ell\geq1}\varGamma^m_{k,\ell}\bigl(\varPhi_{\ell,t}\varPhi_{k,s}-\varPhi_{\ell,s}\varPhi_{k,t}\bigr),
\end{equation}
where the last equality follows by pairing the terms indexed by $(k,\ell)$ and $(\ell,k)$. 

By \eqref{eq:valPhi}, each of the two products in the minor $\varPhi_{\ell,t}\varPhi_{k,s}-\varPhi_{\ell,s}\varPhi_{k,t}$ has valuation at least $2s+2t-k-\ell$. Hence, by the inductive hypothesis, each term on the right-hand side of \eqref{eq:Gammarec} has valuation at least
\[
(3m+10+2k+2\ell)+(2s+2t-k-\ell)=3m+10+2s+2t+k+\ell\geq3(m+1)+10+2s+2t,
\]
since $k+\ell\geq3$. 

This completes the induction.
\end{proof}

\begin{proof}[Proof of Theorem~\ref{thm:internal}]
We fix $m\geq2$ and put $e_m:=3\lfloor(m-1)/2\rfloor+7$. By Lemma~\ref{lem:refined}, we have
\[
\nu(x_{m,1})=e_m\qquad\text{and}\qquad\nu(x_{m+2,1})=e_m+3,
\]
so the normalized first coordinates $x_{m,1}/2^{e_m}$ and $x_{m+2,1}/2^{e_m+3}$ are odd integers. We may therefore choose an odd integer $\lambda_{m+2}$ such that
\begin{equation}\label{eq:lambda}
\frac{x_{m+2,1}}{2^{e_m+3}}\equiv\lambda_{m+2}\,\frac{x_{m,1}}{2^{e_m}}\pmod{2^{3m+16}}.
\end{equation}

By \eqref{eq:valGamma}, we have $\nu(\varGamma^m_{s,1})\geq3m+16$ for any $s\geq2$, and $\varGamma^m_{1,1}=0$. Dividing \eqref{eq:Gamma} with $t=1$ by $2^{e_m}$, we obtain that for any $s\geq1$,
\[
x_{m+2,s}\,\frac{x_{m,1}}{2^{e_m}}\equiv8\,x_{m,s}\,\frac{x_{m+2,1}}{2^{e_m+3}}\pmod{2^{3m+16-e_m}}.
\]
Substituting \eqref{eq:lambda} and cancelling the unit $x_{m,1}/2^{e_m}$ yields
\begin{equation}\label{eq:coeffcong}
x_{m+2,s}\equiv8\lambda_{m+2}\,x_{m,s}\pmod{2^{3\lfloor m/2\rfloor+12}},
\end{equation}
since $3m+16-e_m=3\lfloor m/2\rfloor+12$ for either parity of $m$. 

Finally, substitute \eqref{eq:coeffcong} into \eqref{eq:gfm}. Since $\varrho$ and $\xi$ have integer coefficients, we obtain
\[
\sum_{n\geq0}\PDO\bigl(3\times2^{m+4}n\bigr)q^n\equiv8\lambda_{m+2}\sum_{n\geq0}\PDO\bigl(3\times2^{m+2}n\bigr)q^n\pmod{2^{3\lfloor m/2\rfloor+12}}.
\]
Writing $k=m+2\geq4$, so that $3\lfloor m/2\rfloor+12=3\lfloor k/2\rfloor+9$, proves \eqref{eq:internal}.
\end{proof}

\section{Concluding remarks}\label{sec:final}

For any fixed $m$, the largest $M$ for which $\bx_{m+2}$ is congruent to a scalar multiple of $\bx_m$ modulo $2^M$ equals $\min_{s\geq2}\nu(\varGamma^m_{s,1})-\nu(x_{m,1})$. For $1\leq m\leq9$, we find \[M=14,17,21,24,28,31,35,38,42,\] respectively. Arguing as in the proof of Theorem~\ref{thm:internal}, we conclude that for $3\leq k\leq11$, there exists an odd integer $\lambda_k$ such that, for any $n\geq0$,
\begin{equation}\label{eq:sharp}
\PDO\bigl(3\times2^{k+2}n\bigr)\equiv8\lambda_k\,\PDO\bigl(3\times2^{k}n\bigr)\pmod{2^{\lfloor7(k+1)/2\rfloor}},
\end{equation}
and that the modulus in \eqref{eq:sharp} is best possible for these values of $k$.

This motivates us to make the following conjecture.
\begin{conjecture}\label{conj:sharp}
For any $k\geq3$, there exists an odd integer $\lambda_k$ such that \eqref{eq:sharp} holds for any $n\geq0$.
\end{conjecture}

\noindent It seems to us that the method of Section~\ref{sec:proof3} cannot reach \eqref{eq:sharp} without further ideas.

\medskip

\subsection*{AI Disclosure} The statements of Theorems \ref{thm:parity}, \ref{thm:mod3class} and  \ref{thm:internal}, and the conjecture in Section \ref{sec:final} were found using Anthropic's Claude Fable 5.2. The author already had the auxiliary results and most of the proofs in Sections \ref{sec:preliminaries}--\ref{sec:proof12}, when \cite{LT} was posted on arXiv, which contained \cite{T2} as a reference. This prompted the author to contact Tang who confirmed that he had proved \eqref{eq:strong}. The current paper was then rewritten to emphasize this and the author looked for further congruences and generalizations.

We describe below where AI was used and how\footnote{This is based on a recent template by the `AI+TCS Working Group' held at the Simons Institute: \url{https://simons.berkeley.edu/media/31229/download}.}.

\begin{enumerate}
\item \emph{Asking the question:} The initial scientific questions and prompting for the search space were all contributed by the author, and was inspired by several similar works of the author as well as other mathematicians.

\item \emph{Coming up with the approach and proof development:} The proof techniques used in Sections \ref{sec:preliminaries}--\ref{sec:proof12} are either standard in this area or based on a recently co-authored work of the author \cite{DST}. Fable 5.2 suggested the approach for the proofs of Theorems \ref{thm:parity} and \ref{thm:mod3class} when it conjectured them (prompted by the author to look at other sub-progressions to better the modulus of \eqref{eq:strong}). Once Theorems \ref{thm:parity} and \ref{thm:mod3class} were proved, it was easy to see how to prove the corollaries. The proof strategy for Theorem \ref{thm:internal} was suggested by the author, and then implemented by Fable 5.2.

\item \emph{Writing and exposition:} Some of the initial proofs in Sections \ref{sec:progressions} and the complete initial Section \ref{sec:proof3} were written by Fable 5.2 and then edited by the author for clarity.

\item \emph{Final review:} We used Anthropic's Claude Opus 5.5 to review the final version of this paper. It caught an easily fixable error in the proof of Lemma \ref{lem:Uzeta} and suggested some other minor changes.
\end{enumerate}

Needless to say, the author takes full responsibility for the mathematical content of this paper.

\subsection*{Acknowledgments} We thank Professor Dazhao Tang for sending us his preprint \cite{T2}, and for his helpful comments on an earlier version of this paper.

\subsection*{Declarations} 
\begin{enumerate}
    \item The author declare no conflicts of interest.
    \item The author has no relevant financial or non-financial interests to disclose.
    \item There are no datasets associated with this study.
\end{enumerate}

\bibliographystyle{alpha}

\newcommand{\etalchar}[1]{$^{#1}$}

\end{document}